\documentclass[12pt]{amsart}
\usepackage{graphicx,verbatim,amssymb}
\usepackage[active]{srcltx}

\newtheorem*{thmM}{Main Theorem}

\theoremstyle{definition}

\theoremstyle{remark}

\numberwithin{equation}{section}

\newcommand{\R}{\mathbb R}

\def\Vc{\mathcal{V}}

\def\0{\varnothing}
\def\sm{\setminus}
\def\ol{\overline}
\def\wt{\widetilde}
\def\Vol{{\rm Vol}}

\renewcommand\le{\leqslant}
\renewcommand\ge{\geqslant}

\begin{document}

\title[Co-convex Aleksandrov--Fenchel inequality]
{Aleksandrov--Fenchel inequality for coconvex bodies}
\author{Askold Khovanski\u{\i}}
\author{Vladlen Timorin}

\address[Askold Khovanskii]
{Department of Mathematics, University of Toronto, Toronto,
Canada; Moscow Independent University; Institute for Systems Analysis,
Russian Academy of Sciences.
}

\email{askold@math.toronto.edu}

\address[Vladlen Timorin]{Faculty of Mathematics and Laboratory of Algebraic Geometry,
National Research University Higher School of Economics,
7 Vavilova St 117312 Moscow, Russia; 
Independent University of Moscow,
Bolshoy Vlasyevskiy Pereulok 11, 119002 Moscow, Russia}

\email{vtimorin@hse.ru}

\thanks{
The first named author was partially supported by Canadian Grant \textrm{N~156833-02}.
The second named author was partially supported by
the Dynasty Foundation grant, RFBR grants 11-01-00654-a, 12-01-33020,
and AG Laboratory NRU HSE, MESRF grant ag. 11 11.G34.31.0023
}


\begin{abstract}
We prove a version of the Aleksandrov--Fenchel inequality for
mixed volumes of coconvex bodies.
This version is motivated by an inequality from commutative algebra
relating intersection multiplicities of ideals.
\end{abstract}
\maketitle

\section{Introduction}
We start by recalling the classical Aleksandrov--Fenchel inequality.
The \emph{Minkowski sum} of two convex sets $A$, $B\subset \R^d$ is defined as
$A+B=\{a+b\,|\, a\in A,\ b\in B\}$.
For a positive real number $\lambda$, we let $\lambda A$ denote the
set $\{\lambda a\,|\, a\in A\}$.
By definition, a \emph{convex body} is a compact convex set, whose interior
is nonempty.
A \emph{linear family of convex bodies} is a collection of the following
objects: a real vector space $\Vc$, an open subset $\Omega\subset\Vc$,
a map $f$ from $\Omega$ to the set of all convex bodies in $\R^d$ such that
$$
f(\lambda_1 v_1+\dots+\lambda_n v_n)=\lambda_1 f(v_1)+\dots +\lambda_n f(v_n)
$$
whenever all $v_i\in\Omega$,
all $\lambda_i$ are positive, and $\lambda_1 v_1+\cdots+ \lambda_n v_n\in\Omega$.
A linear family of convex bodies \emph{with $m$ marked points} is a linear family $(\Vc,\Omega,f)$ of convex bodies, in which some $m$ elements of $\Omega$ are marked.

With every linear family $\alpha=(\Vc,\Omega,f)$ of convex bodies, we associate
the \emph{volume polynomial} $\Vol_\alpha$ as follows.
For $v\in\Omega$, we define $\Vol_\alpha(v)$ as the usual $d$-dimensional
volume of the convex body $f(v)$.
It is well known that the function $\Vol_\alpha$ thus defined extends
to a unique polynomial on $\Vc$ that is homogeneous of degree $d$.
For $v\in\Vc$, we let $L_v$ denote the usual directional (Lie) derivative along $v$.
Thus $L_v$ is a differential operator that acts on functions, in particular,
degree $k$ polynomials on $\Vc$ are mapped by this operator to degree $k-1$
polynomials.
If $\alpha=(\Vc,\Omega,f)$ is a linear family
of convex bodies with $(d-2)$ marked points $v_1$, $\dots$, $v_{d-2}\in\Omega$,
then we define the \emph{Aleksandrov--Fenchel} symmetric bilinear form $B_\alpha$ on $\Vc$
by the formula
$$
B_\alpha(u_1,u_2)=\frac 1{d!}L_{u_1}L_{u_2}L_{v_1}\dots L_{v_{d-2}}(\Vol_\alpha).
$$
Note that the expression in the right-hand side is a real number
(called the \emph{mixed volume} of the convex bodies $f(u_1)$, $f(u_2)$, $f(v_1)$, $\dots$,
$f(v_{d-2})$).
Indeed, this is the image of a homogeneous degree $d$ differential operator
with constant coefficients on a homogeneous degree $d$ polynomial.
The corresponding quadratic form $Q_\alpha(u)=B_\alpha(u,u)$ is given by
the formula $Q_\alpha=\frac 2{d!}L_{v_1}\dots L_{v_{d-2}}(\Vol_\alpha)$.
The expression in the right-hand side is the result of the action of
a degree $d-2$ homogeneous differential operator with constant coefficients on a homogeneous
degree $d$ polynomial, i.e. a quadratic form.
\emph{The Aleksandrov--Fenchel inequality} \cite{Al} states that, for all $u_1\in\Vc$ and
$u_2\in\Omega$, we have
$$
B_\alpha(u_1,u_2)^2\ge B_\alpha(u_1,u_1) B_\alpha(u_2,u_2).
$$
The Aleksandrov--Fenchel inequality is a far-reaching generalization of
the classical isoperimetric inequality. See \cite{Mc,VT} for generalizations
of the Aleksandrov--Fenchel inequality for convex polytopes.

Let $C\subset\R^d$ be a strictly convex cone with the apex at $0$
and a nonempty interior.
A subset $A\subset C$ is called ($C$-)\emph{coconvex}
if the complement $A\sm C$ is convex.
Define a ($C$-)\emph{coconvex body} as a compact
coconvex subset of $C$ with a nonempty interior.
If $A$ and $B$ are coconvex sets with respect to the same cone $C$,
then we can define $A\oplus B$ as the complement in $C$ of
the Minkowski sum of $C\sm A$ and $C\sm B$.
Define a \emph{linear family of $C$-coconvex bodies} as a collection of
the following objects: a vector space $\Vc$,
an open subset $\Omega\subset\Vc$, a map $g$ from $\Omega$ to the set
of all $C$-coconvex bodies such that
$$
g(\lambda_1 v_1+\dots +\lambda_n v_n)=\lambda_1 g(v_1)\oplus\dots \oplus\lambda_n g(v_n)
$$
whenever all $v_i\in\Omega$, all $\lambda_i$ are positive, and
$\lambda_1 v_1+\dots \lambda_n v_n\in\Omega$.
A linear family of $C$-coconvex bodies \emph{with $m$ marked points} is a 
linear family $(\Vc,\Omega,g)$ of $C$-coconvex bodies, 
in which some $m$ elements of $\Omega$ are marked.
With every linear family $\beta$ of $C$-coconvex bodies, we associate the
volume function $\Vol_\beta$ in the same way as with a linear family of convex bodies.
The function $\Vol_\beta$ thus defined is also a homogeneous degree $d$ polynomial
(we will prove this below).
Given a linear family $\beta$ of $C$-coconvex bodies with $d-2$ marked points
$v_1$, $\dots$, $v_{d-2}$, we define the \emph{coconvex Aleksandrov--Fenchel}
symmetric bilinear form as
$$
B_\beta^C(u_1,u_2)=\frac 1{d!}L_{u_1}L_{u_2}L_{v_1}\dots L_{v_{d-2}}(\Vol_\beta).
$$
We will also consider the corresponding quadratic form
$Q^C_\beta=\frac 2{d!}L_{v_1}\dots L_{v_{d-2}}(\Vol_\beta)$.
Our main result is the following

\begin{thmM}
  \label{t:main}
  The form $Q_\beta^C$ is non-negative, i.e. $Q_\beta^C(u)\ge 0$
  for all $u\in\Vc$.
  In particular, the corresponding symmetric bilinear form
  satisfies the Cauchy--Schwartz inequality
  $$
  B_\beta^C(u_1,u_2)^2\le B_\beta^C(u_1,u_1) B_\beta^C(u_2,u_2).
  $$
\end{thmM}

The inequality stated in the Main Theorem is called the \emph{coconvex
Aleksandrov--Fenchel inequality}.
In recent paper \cite{Fi}, Theorem \ref{t:main} is proved under the assumption
that $C$ is a fundamental cone of some Fuchsian group $\Gamma$ acting by linear isometries
of a pseudo-Euclidean metric, and $C\sm A_k$ is the intersection of some convex
$\Gamma$-invariant set with $C$.

Theorem \ref{t:main} is motivated by an Aleksandrov--Fenchel type inequality
for (mixed) intersection multiplicities of ideals \cite{KK}.
There was an earlier preprint of ours, in which more sophisticated tools were used
for proving Theorem \ref{t:main}. We plan a sequel to this paper that will deal
with representations of coconvex bodies as virtual convex bodies in the sense of \cite{KP}.

The following inequalities follow from Theorem \ref{t:main} in the same way as
similar inequalities for convex bodies follow from the classical Alexandrov--Fenchel
inequality (cf. \cite{Fi}):

\noindent\emph{Reversed Brunn--Minkowski inequality:} the function $\Vol_\beta^{\frac 1d}$ is convex, i.e.
$$
(\Vol_\beta(tu+(1-t)v))^\frac 1d\le t\Vol_\beta(u)^\frac 1d+(1-t)\Vol_\beta(v)^\frac 1d,\quad t\in [0,1].
$$

\noindent\emph{Generalized reversed Brunn--Minkowski inequality:} the function
$\left(L_{v_1}\dots L_{v_k}\Vol_\beta\right)^\frac 1{d-k}$ is convex.

\noindent\emph{First reversed Minkowski inequality:}
$$
\left(\frac 1{d!}L_{u}L_{v}^{d-1}(\Vol_\beta)\right)^d\le \Vol_\beta(u)\Vol_\beta(v)^{d-1}.
$$

\noindent\emph{Second reversed Minkowski inequality:} if all marked points coincide with $u$, then
$$
B_\beta^C(u,v)^2\le \Vol_\beta(u)\, B_\beta^C(v,v)
$$

\section{Proof of the Main Theorem}
Recall that every quadratic form $Q$ on a finite dimensional real vector space can be represented in
the form
$$
x_1^2+\cdots+x_k^2-x_{k+1}^2-\cdots-x_{k+\ell}^2
$$
for a suitable linear coordinate system $(x_1,\dots,x_m)$, $m\ge k+\ell$.
The pair $(k,\ell)$ is called the \emph{signature} of $Q$.
It is well known that $Q$ has signature $(1,\ell)$ for some $\ell$ if and only if
\begin{enumerate}
  \item there exists a vector $v_0\in\Vc$ with $Q(v_0)>0$;
  \item the corresponding symmetric bilinear form $B$ (such that $B(u,u)=Q(u)$)
  satisfies the reversed Cauchy--Schwartz inequality: $B(u,v)^2\ge Q(u)Q(v)$
  for all $u\in\Vc$ and all $v\in\Vc$ such that $Q(v)>0$.
\end{enumerate}
Thus, the Aleksandrov--Fenchel inequality is equivalent to the fact that
$Q_\alpha$ has signature $(1,\ell(\alpha))$ for every finite-dimensional linear
family $\alpha$ of convex $d$-dimensional bodies with $d-2$ marked points.

Fix a strictly convex cone $C$ with the apex at $0$ and a $C$-linear family
$\beta=(\Vc,\Omega,g)$ of coconvex bodies.
We may assume that $g(\Omega)$ is bounded in the sense that there is
a large ball in $\R^d$ centered at $0$ that contains all
coconvex bodies $g(v)$, $v\in\Omega$.
Then there is a linear functional $\xi$ on $\R^d$ and a number $t_0>0$
such that $g(v)$ lies in the half-space
$$
W(t_0)=\{x\in\R^d\,|\,\xi(x)\le t_0\},
$$
for every $v\in\Omega$.
Choose any $t_1>t_0$.
We will now define a linear family
$\alpha=(\Vc\times\R,\Omega\times (t_0,t_1),f)$ of convex bodies as follows.
For $v\in\Omega$ and $t\in (t_0,t_1)$, we set $f(v,t)$ to be the convex
body $\ol{(C\cap W(t))\sm g(v)}$.
The proof of the coconvex Aleksandrov--Fenchel inequalities is based
on the comparison between the linear families $\alpha$ and $\beta$.

Let $\pi:\Vc\times\R\to\Vc$ denote the natural projection.
The polynomials $\Vol_\alpha$ and $\Vol_\beta\circ\pi$ are defined on the same space
$\Vc\times\R$.
We have the following relation between the two polynomials:
$$
\Vol_\alpha=\Vol(C\cap W(t))-\Vol_\beta\circ\pi,\eqno{(V)}
$$
which is clear from the additivity of the volume.
The first term in the right-hand side has the form $ct^d$, where $c$
is some positive constant.
The second term in the right-hand side does not depend on $t$.
It follows from $(V)$ that $\Vol_\beta$ is a homogeneous degree $d$ polynomial.

Let us mark some points $(v_1,s_1)$, $\dots$, $(v_{d-2},s_{d-2})$ in
$\Omega\times (t_0,t_1)$.
Apply the differential operator $\frac 2{d!}L_{(v_1,s_1)}\dots L_{(v_{d-2},s_{d-2})}$
to both sides of $(V)$.
We obtain that
$$
Q_\alpha=c't^2-Q^C_\beta\circ\pi,\eqno{(Q)}
$$
where $c'$ is some positive constant (equal to $c s_1\cdots s_{d-2}$).
The last term in the right-hand side of $(Q)$ is obtained using the chain rule
and the fact that the differential of $\pi$ coincides with $\pi$.
Consider the following linear family $\wt\beta=(\Vc\times\R,\Omega\times (t_0,t_1),\wt g)$
of $C$-coconvex polytopes: $\wt g=g\circ\pi$.
Then $Q^C_{\wt\beta}=Q^C_\beta\circ\pi$, again by the chain rule.
In the right-hand side of $(Q)$, we have the difference of two quadratic
forms, moreover, these two forms depend on disjoint sets of variables.

If $q_1$, $q_2$ are quadratic forms depending on disjoint sets of variables,
and $(k_1,\ell_1)$, $(k_2,\ell_2)$, respectively, are signatures of these forms,
then $q_1+q_2$ is a quadratic form of signature $(k_1+k_2,\ell_1+\ell_2)$.
We now apply this observation to identity $(Q)$.
The first term of the right-hand side, $c't^2$, has signature $(1,0)$.
The signature of the left-hand side is equal to $(1,\ell)$ for some $\ell\ge 0$,
by the classical Alexandrov--Fenchel inequality.
It follows that the signature of $Q^C_{\wt\beta}$ is $(\ell,0)$, i.e. the
form $Q^C_{\wt\beta}$ is non-negative.
Finally, since $Q^C_{\wt\beta}=Q^C_\beta\circ\pi$, the quadratic form
$Q^C_\beta$ is also non-negative.

\end{document}